\documentclass[12pt]{amsart}

\newtheorem{THM}{Theorem}[section]
\newtheorem{LMA}[THM]{Lemma}
\newtheorem{PROP}[THM]{Proposition}
\newtheorem{CORO}[THM]{Corollary}
\newtheorem{CONJ}[THM]{Conjecture}
\newtheorem{EG}[THM]{Example}

\numberwithin{equation}{section}

\usepackage{amsmath}
\usepackage{amssymb}
\usepackage{euscript}

\newcommand{\showon}{\begin{eqnarray*}}
\newcommand{\showoff}{\end{eqnarray*}}

\newcommand{\none}{\varnothing}
\newcommand{\drop}{\smallsetminus}
\newcommand{\diag}{\mathrm{diag}}
\newcommand{\goesto}{\rightarrow}
\newcommand{\one}{\boldsymbol{1}}
\newcommand{\zero}{\boldsymbol{0}}
\newcommand{\Cov}{\mathrm{Cov}}
\newcommand{\Supp}{\mathrm{Supp}}
\newcommand{\la}{\langle}
\newcommand{\ra}{\rangle}

\newcommand{\tri}{\triangle}
\newcommand{\rank}{\mathrm{rank}}

\newcommand{\A}{\EuScript{A}} 
\newcommand{\B}{\EuScript{B}}
\newcommand{\BB}{\mathbf{B}}

\newcommand{\G}{\EuScript{G}}

\newcommand{\II}{\mathbf{I}}

\renewcommand{\L}{\EuScript{L}}
\newcommand{\M}{\EuScript{M}} \newcommand{\m}{\mathbf{m}}
\newcommand{\N}{\EuScript{N}}  \newcommand{\NN}{\mathbb{N}}

\newcommand{\Q}{\EuScript{Q}} 

\newcommand{\R}{\EuScript{R}} 
 \renewcommand{\S}{\EuScript{S}}
\renewcommand{\SS}{\mathbf{S}}

\newcommand{\U}{\EuScript{U}}

\newcommand{\w}{\mathbf{w}}

 \newcommand{\y}{\mathbf{y}}

\newcommand{\z}{\mathbf{z}} \newcommand{\zZ}{\widetilde{Z}}

\begin{document}

\title[Negative correlation]{Negatively correlated random variables
and Mason's Conjecture}

\author{David G. Wagner}
\address{Department of Combinatorics and Optimization\\
University of Waterloo\\
Waterloo, Ontario, Canada\ \ N2L 3G1}
\email{\texttt{dgwagner@math.uwaterloo.ca}}
\thanks{Research supported by the Natural
Sciences and Engineering Research Council of Canada under
operating grant OGP0105392.}
\keywords{matroid, delta--matroid, logarithmic concavity,
Rayleigh monotonicity, Potts model.}
\subjclass{05A20;\ 05B35, 60C05, 82B20.}

\date{}

\begin{abstract}
Mason's Conjecture asserts that for an $m$--element rank $r$ matroid $\M$
the sequence $(I_k/\binom{m}{k}:\ 0\leq k\leq r)$ is logarithmically concave,
in which $I_k$ is the number of independent $k$--sets of $\M$.  A related
conjecture in probability theory implies these inequalities provided that
the set of independent sets of $\M$ satisfies a strong negative correlation
property we call the \emph{Rayleigh condition}.  This condition is known to
hold for the set of bases of a regular matroid.  We show that if $\omega$ is
a weight function on a set system $\Q$ that satisfies the Rayleigh condition
then $\Q$ is a convex delta--matroid and $\omega$ is logarithmically submodular.
Thus, the hypothesis of the probabilistic conjecture leads inevitably to matroid
theory.  We also show that two--sums of matroids preserve the Rayleigh condition
in four distinct senses, and hence that the Potts model of an iterated two--sum of
uniform matroids satisfies the Rayleigh condition.  Numerous conjectures and auxiliary
results are included.
\end{abstract}

\maketitle

\section{Introduction.}

Mason's Conjecture \cite{Ma} is that the sequence $(I_k:\ 0\leq k\leq r)$ of numbers
of independent $k$--sets of an $m$--element rank $r$ matroid $\M$ is
logarithmically concave in the strong sense (I-4) that
$(I_k/\binom{m}{k}:\ 0\leq k\leq r)$ is log--concave.  That is, that
$$\frac{I_{k}^{2}}{\binom{m}{k}^{2}}\geq
\frac{I_{k-1}}{\binom{m}{k-1}}\cdot\frac{I_{k+1}}{\binom{m}{k+1}}$$
for all $1\leq k\leq r-1$.  A weaker form of the conjecture is (I-2) that
the sequence $(I_k:\ 0\leq k\leq r)$ itself is log--concave:\
$I_k^2\geq I_{k-1} I_{k+1}$ for all $1\leq k\leq r-1$.
Mahoney \cite{Mah} has shown that (I-2) holds for graphic (cycle) matroids of
outerplanar graphs.
Dowling \cite{Dowl} has shown the inequalities $I_k^2\geq I_{k-1}I_{k+1}$
in general for $1\leq k\leq 7$. 
Zhao \cite{Zhao} has shown that $I_k^2\geq (1+1/k)I_{k-1}I_{k+1}$
in general for $1\leq k\leq 5$.
These are currently the most notable partial results on Mason's Conjecture.

There is a related conjecture in probability theory, but its origin is
obscure.  Pemantle \cite{Pe} considers a lot of conditions of this kind.
The \emph{Big Conjecture} 3.4 states that if $\omega:\B(E)\goesto [0,\infty)$ is a
nonnegative weight function on a finite Boolean algebra $\B(E)$, and if $f_k(\omega):=
\sum_{S\subseteq E:\ |S|=k}\omega(S)$
for all $0\leq k\leq m=|E|$ then $(f_k(\omega)/\binom{m}{k}:\ 0\leq k\leq m)$ is
logarithmically concave, provided that $\omega$ satisfies something we call the
\emph{Rayleigh condition}.  This condition is a strong pairwise negative correlation
property among random variables $\{X_e:\ e\in E\}$ corresponding to the elements
of the ground--set $E$, with joint distribution function encoding the weight function
$\omega$.  The Rayleigh condition is known to hold in its weakest form
($\BB$--Rayleigh, for bases) for all regular (unimodular) matroids, and for many more
\cite{CW}.  There are more refined and informative versions of the Rayleigh
condition for matroids:\ $\II$--Rayleigh, $\SS$--Rayleigh, and Potts--Rayleigh
for independent sets, spanning sets, and the Potts model, respectively. 

A positive solution to the Big Conjecture would be a very good thing.
If so, then every $\II$--Rayleigh matroid satisfies Mason's Conjecture (I-4).
In Section 5 we see that every series--parallel matroid is $\II$--Rayleigh, and
we have reason to believe that the class of $\II$--Rayleigh matroids 
might contain all graphs, maybe all regular matroids, perhaps even
more.  Thus, this line of reasoning has the potential for substantial progress
on Mason's Conjecture.  

Although the Big Conjecture has not been proven
we do have a new equivalent form of it, Conjecture 3.11,
which states that if $\omega$
satisfies the Rayleigh condition then its \emph{symmetrization}
$\widetilde{\omega}$ also satisfies the Rayleigh condition. 
By the exchangeable (symmetric function) case of the Big Conjecture
-- that is, Proposition 3.6 -- this implies the inequalities on $(f_k(\omega))$.
This suggests an entirely different approach towards the required inequalities.

In Section 2 we briefly review some unimodality conditions for nonnegative
real sequences, some sequences associated with matroids, and some relevant
unimodality conjectures and results.  This is meant to put the results of
later sections in context.

In Section 3 we look at some examples,
state the Big Conjecture 3.4,
prove the \emph{exchangeable case} Proposition 3.6 of it, and
review some supplementary results.
This is also partly a capsule summary of some of Section 2.4 of Pemantle \cite{Pe}.
Then we give a new equivalent form of the Big Conjecture 3.11,
and, after some algebra, the sufficient conditions Conjectures 3.13 and 3.14.
These latter two conjectures are more local than 3.11, so even though they are
strictly stronger they might be more amenable to proof.

In Section 4 we show that if $\omega$ is Rayleigh then $\Supp(\omega)$,
the set of sets on which $\omega$ is positive, is a convex delta--matroid,
and that $\omega$ is logarithmically submodular.
Regarding the conjectures of Section 2 this is a negative result:\
the Big Conjecture is directly relevant only to Mason's Conjecture (I-4).
On the other hand, this structure might be useful in an attempt to prove
the Big Conjecture.

In Section 5 we turn to finding examples to populate the theory.  We see 
that uniform matroids are Potts--Rayleigh.  We show that the Rayleigh
condition on the Potts model partition function
is preserved by two--sums of matroids.  Consequently, every 
series--parallel matroid is Potts--Rayleigh.
Analogously, two--sums preserve
the Rayleigh condition for matroids in any of the three frozen senses:\
for bases, for independent sets, or for spanning sets.  Concerning the
$\II$--Rayleigh property for graphs we have a few 
small examples and two relatively technical conjectures. 
CJSSS \cite{CJSSS} gives the generating function
for the set of spanning forests of a graph as a Grassmann--Berezin integral.
This is a beautiful result, and can only help any attempt to prove that
graphs are $\II$--Rayleigh.
As an adjunct to the Big Conjecture 3.4 we give a related scale of Conjectures 5.11,
guessing that various classes of matroids are $\II$--Rayleigh.
Some binary matroids fail to be balanced \cite{SW}, but within the class of
sixth--root--of--unity matroids there are no show--stoppers in sight.
They might all be Potts--Rayleigh!  The relationship between the
Potts--Rayleigh condition and the half--plane property (HPP) is unclear --
my guess is for counterexamples both ways.
There is still very little data to go on, and many interesting examples are
sure to be as yet undiscovered.

In preparing this paper I have benefitted from conversations and correspondence
with many people.  In particular I thank
Andr\'e Bouchet,
Seth Chaiken,
Bill Cunningham,
Jim Geelen,
Bill Jackson,
Tom Liggett,
Robin Pemantle,
and Alan Sokal
for their comments.
Thanks also to Marc Noy for organizing a very successful 
\emph{2nd Workshop on Tutte Polynomials} at the
Universitat Aut\`onoma de Barcelona, Oct. 4--7, 2005, at which an
early version of this paper was presented.

\section{Logarithmic Concavity Conjectures for Matroids.}

2.1. \ \textsc{unimodality conditions.}\\

Let $0\leq s\leq r\leq m$ be integers, and let $a_{s}, a_{s+1},\ldots, a_{r}$ be a
finite sequence of nonnegative real numbers.  Consider the following conditions on 
this sequence $a=(a_{k}:\ s\leq k\leq r)$:\\
(a-0)\ \emph{no internal zeros}:\ if $s\leq i<j<k\leq r$ and 
$a_{i}a_{k}\neq 0$, then $a_{j}\neq 0$;
(a-1)\  \emph{unimodality}:\ $a_{s}\leq a_{s+1}\leq \cdots \leq a_{p} 
\geq \cdots \geq a_{r}$ for some $s\leq p\leq r$;\\
(a-2)\ \emph{logarithmic concavity}:\ 
$a_{k}^{2}\geq a_{k-1}a_{k+1}$ for all $s+1\leq k\leq r-1$;\\
(a-3)\  logarithmic concavity of the sequence $(k!a_{k}:\ s\leq k\leq 
r)$;\\
(a-4)\  logarithmic concavity of the sequence $(a_{k}/\binom{m}{k}:\ 
s\leq k\leq r)$;\\
(a-5)\ logarithmic concavity of the sequence $(a_{k}/\binom{r}{k}:\ 
s\leq k\leq r)$;\\
(a-6)\  the polynomial $a_{s}+a_{s+1}t+\cdots+ a_{r}t^{r}$ has only real 
(nonpositive) zeros.\\
Elementary arguments show that (a-5) $\Longrightarrow$ (a-4)
$\Longrightarrow$ (a-3) $\Longrightarrow$ (a-2),
and that (a-2) and (a-0) together imply (a-1).  Newton's Inequalities
(item (51) of \cite{HLP}) assert that (a-6) implies both (a-5) and (a-0).
The sequences we consider are usually easily seen to satisfy (0).\\

2.2. \ \textsc{sequences associated with matroids.}\\

Let $\M$ be a loopless matroid of rank $r$ on a set $E$ of size $|E|=m$.
Several sequences associated with $\M$ have been 
conjectured to satisfy one or another of the above conditions.  In 
most cases counterexamples to condition (5) can be found easily,
so condition (4) is the strongest reasonable conjecture.\\

$(W_{k}:\ 0\leq k\leq r)$ in which $W_{k}$ is the number of 
flats of $\M$ of rank $k$.
Unimodality (W-1) was conjectured by Rota \cite{Ro} in the late 1960s, and
logarithmic concavity in any of the forms (W-2) to (W-4) was
conjectured by Mason \cite{Ma} in the early 1970s.  
Seymour \cite{Sey} has shown that $W_{2}^{2}\geq W_{1}W_{3}$ for matroids
with at most four points on any line.

$(I_{k}:\ 0\leq k\leq r)$ in which $I_{k}$ is the number of 
independent sets of $\M$ of size $k$.
Unimodality (I-1) was conjectured by Welsh \cite{We} in the late 1960s, and
logarithmic concavity in any of the forms (I-2) to (I-4) was
conjectured by Mason \cite{Ma} in the early 1970s.
Partial results were reviewed in the first paragraph of the Introduction.
If Conjecture 3.4 is true then every series--parallel graph satisfies (I-4),
which would be progress.

$(\chi_{k}:\ 0\leq k\leq r)$ in which $\chi_{k}$
is the number of subsets of $E$ of size $k$ containing no broken circuit of $\M$
(relative to any fixed total order on $E$).  For graphic matroids these
are the coefficients of the chromatic polynomial of the graph.
Unimodality ($\chi$-1) was conjectured by Read \cite{Re} in the late 1960s (for 
graphs), and logarithmic concavity in the form ($\chi$-2) was conjectured by Hoggar 
\cite{Hog} in the early 1970s (also for graphs).  The literature on zeros or
coefficients of chromatic polynomials is extensive --- see \cite{BRW,Jack,Re,Sok}
for starters.

$(h_{k}:\ 0\leq k\leq r)$ in which the integers $h_{k}$ are 
defined by the relation
$$\sum_{k=0}^{r} I_{k} t^{k} = \sum_{k=0}^{r}h_{k}x^{k}(1+x)^{r-k}.$$
The properties (h-1) to (h-4) were conjectured for this sequence by 
Dawson \cite{Daw} in the early 1980s.  Dawson proves that (h-2,0) implies (I-2,0),
and that the sequence $(h_k/\binom{m}{k}: 0\leq k\leq r)$
is nonincreasing.

Fix any $E'\subseteq E$ and consider $(c_{k}:\ 0\leq k\leq r)$
in which $c_{k}$ is the number of bases $B$ of $\M$ such that $|B\cap E'|=k$.
Condition (c-5) was proven by Stanley \cite{St} and condition (c-6) by Godsil
\cite{Go} in 
the early 1980s, both for the class of regular matroids.
Condition (c-6) was proven recently for the larger class of HPP
matroids by Choe and Wagner \cite{CW};\ see also \cite{COSW,Wa2}.
Does the $\BB$--Rayleigh condition imply (c-2)?\\

It must be said that conjectures ($\chi$-1) and (h-1) now seem dubious
in the generality of all matroids.  Even (I-1) seems a little suspect
since Bj\"orner \cite{Bj1,Bj2} has given counterexamples to
(I-1) and (h-1) in the somewhat wider arena of shellable simplicial complexes.
But we are not asking for universal results -- rather, just for the
identification of significant classes of matroids (or related objects)
satisfying more--or--less restrictive versions of these unimodality conditions.\\

2.3. \  \textsc{from integer sequences to polynomials.}\\ 

For a loopless matroid $\M$ on a set $E$ and a set $\m=\{m_e:\ e\in E\}$
of positive integers indexed by $E$, let $\M[\m]$ be the matroid obtained
from $\M$ by replacing each $e\in E$ by $m_e$ elements in parallel.  The
number of $k$--element independent sets of $\M[\m]$ is
$$I_k(\M[\m])=\sum_{S\in\II\M:\ |S|=k}\ \  \prod_{e\in S} m_e,$$
in which $\II\M$ is the simplicial complex of independent sets of $\M$.
Considering the conjectures about $(I_k)$ for all of these matroids at the
same time, we are led to consider properties of the polynomial
$$Z(\II\M;\y):=\sum_{S\in\II\M}\ \ \prod_{e\in S} y_e,$$
in which $\y=\{y_e:\ e\in E\}$ is a set of algebraically independent commuting
indeterminates. The notation $\y^S:=\prod_{e\in S}y_e$ is useful.

\section{Negatively Correlated Random Variables.}

3.1. \ \textsc{partition functions and the rayleigh condition.}\\

For our purposes it suffices to consider finite sets of binary (zero or one) valued
random variables.

Let $E$ be a finite set with $|E|=m$, let $\B(E)$ be the set of all 
subsets of $E$, and let $\omega:\B(E)\goesto[0,\infty)$ be a
nonnegative--valued function on $\B(E)$ that is not identically zero.
Let $\y:=\{y_{e}:\ e\in E\}$
be a set of algebraically independent commuting indeterminates, and consider the
\emph{partition function}
$$Z(\omega;\y):=\sum_{S\subseteq E}\omega(S)\y^{S}.$$
For any choice of positive values $y_{e}>0$ for each $e\in E$, this determines
a probability measure $\mu=\mu_{\y}$ on $\B(E)$ by setting
$$\mu(S):=\frac{\omega(S)\y^{S}}{Z(\omega;\y)}$$
for all $S\subseteq E$.
The \emph{atomic random variables} of this theory are $X_{e}$ for
each $e\in E$, given by
$$X_{e}(S):=\left\{\begin{array}{ll}
1 & \mathrm{if}\ e\in S,\\
0 & \mathrm{if}\ e\not\in S.
\end{array}\right.$$
The \emph{expectation} of a random variable $X$ is
$$\la X \ra :=\sum_{S\subseteq E} X(S)\mu(S).$$
The \emph{covariance} of two random variables $X$ and $Y$ is
$$\Cov(X,Y):=\la XY \ra - \la X \ra \la Y \ra.$$

The hypothesis we put on the weight function $\omega$ is the 
following:\  for any positive choice of parameters $\y>\zero$,
and any distinct $e\neq f$ in $E$, $\Cov(X_{e},X_{f})\leq 0$.
As a short codename for this hypothesis, we will say that
the weight function $\omega$ or the partition function $Z(\omega;\y)$
\emph{satisfies the Rayleigh condition};  even more briefly, we will 
say that $\omega$ or $Z$ \emph{is Rayleigh}.   The reason for this
terminology is by analogy with the Rayleigh monotonicity property of
(linear resistive) electrical networks. That is (with the notation of
Example 3.1) that $\BB\G$ is Rayleigh for a graphic matroid $\G$.
In fact, $\BB\M$ is Rayleigh for a much wider class of matroids than the
regular ones \cite{Choe1,COSW,CW,Wa1,Wa2}.

The following three examples define classes of weight functions $\omega$
for which the question ``Is $\omega$ Rayleigh?'' seems interesting.

\begin{EG}\emph{
Let $\Q\subseteq \B(E)$ be any set system, and define
$\one_\Q:\B(E)\goesto[0,\infty)$ by
$$\one_\Q(S):=\left\{\begin{array}{ll}
1 & \mathrm{if}\ S\in\Q,\\
0 & \mathrm{if}\ S\not\in\Q.
\end{array}\right.$$
The corresponding partition function
$Z(\one_\Q;\y)$ will be denoted more simply by $Z(\Q;\y)$.
In particular, for a matroid $\M$ we consider this construction with
$\Q$ being any of:\\
the set $\BB\M$ of bases of $\M$,\\
the set $\II\M$ of independent sets of $\M$,\\
the set $\SS\M$ of spanning sets of $\M$.\\
}\end{EG}

\begin{EG}\emph{
For a matroid $\M$ of rank $r$ on the ground set $E$, and real $q>0$, define
$\omega:\B(E)\goesto(0,\infty)$ by
$$\omega(S):=q^{-\mathrm{rank}_{\M}(S)}.$$
Denote the corresponding partition function by $Z(\M,q;\y)$.
It is the partition function of the \emph{$q$--state Potts model} associated with
the matroid $\M$.  The set systems of Example $3.1$ are limiting\
cases of this example, as follows.  Fix $0\leq\alpha\leq 1$ and consider
the substitution
$q^{(1-\alpha)r}Z(\M,q;q^\alpha\y)$.  For $S\subseteq E(\M)$, the
coefficient of $\y^S$ in this is $q$ to the power
$$(1-\alpha)(r-\mathrm{rank}_\M(S))+\alpha(|S|-\mathrm{rank}_\M(S)).$$
This exponent is nonnegative for all $S\subseteq E(\M)$.  As $q\goesto 0$
only those terms in which the exponent is zero survive.  Thus
$$\lim_{q\goesto 0}q^{(1-\alpha)r}Z(\M,q;q^{\alpha}\y)=
\left\{\begin{array}{ll}
Z(\SS\M;\y) & \mathrm{if}\ \alpha=0,\\
Z(\BB\M;\y) & \mathrm{if}\ 0<\alpha<1,\\
Z(\II\M;\y) & \mathrm{if}\ \alpha=1.
\end{array}\right.$$
See Sokal \cite{Sok} for an excellent survey of the combinatorial and
analytic properties of Potts model partition functions.  Limiting
arguments can be used to show that if $Z(\M,q;\y)$ is Rayleigh for all
$q$ in an interval $0<q\leq q_*(\M)$ then both $\II\M$ and $\SS\M$ are
Rayleigh, and that if either $\II\M$ or $\SS\M$ is Rayleigh then
$\BB\M$ is Rayleigh.
}\end{EG}

\begin{EG}\emph{
A \emph{nonsingular M--matrix} is a square real symmetric matrix for which
all principal minors are positive and all off--diagonal elements are
nonpositive.  If either $A$ or $A^{-1}$ is a nonsingular M--matrix with rows
and columns indexed by $E$, define $\omega_A:\B(E)\goesto(0,\infty)$ by putting
$\omega_A(S)$ equal to
the minor of $A$ indexed by rows and columns in $S$, for every $S\subseteq E$.
(By Jacobi's formula for complementary minors of inverse matrices, this is possible.)
With the notation of Conjecture 3.4, Holtz \cite{Hol} has recently proven
(f-4,0) for such $\omega_A$.
}\end{EG}

\vspace{4mm}

3.2. \  \textsc{the big conjecture.}\\

\begin{CONJ}
Let $Z(\omega;\y)$ satisfy the Rayleigh condition.  For each $0\leq 
k\leq m$ let
$$f_{k}(\omega):=\sum_{S\subseteq E:\ |S|=k}\omega(S).$$
Then $(f_{k}(\omega)/\binom{m}{k}:\ 0\leq k\leq m)$ is 
logarithmically concave with no internal zeros, the condition
\textup{(f-4,0)}.
\end{CONJ}

Note that the $\y$ can be included in the conclusion by considering
$f_k(\omega;\y):=\sum_{S\subseteq E:\ |S|=k}\omega(S)\y^S$.
But the hypothesis is unchanged by rescaling $y_e=w_e z_e$ for each
$e\in E$. Considering $Z(\omega; \w\z)$ as a polynomial in $\z$,
the Big Conjecture implies (f$(\w)$-4,0) for all $\w>\zero$.
Thus, nothing is gained.

Note also that Proposition 4.7(a) implies that if $Z(\omega;\y)$ is Rayleigh
then $(f_{k}(\omega):\ 0\leq k\leq m)$ has no internal zeros,
the condition (f-0).  Since Proposition 4.7 is derived independently of 
the results of this section we may make use of it here.

If Conjecture 3.4 is true then any matroid $\M$ for which
$\II\M$ is Rayleigh satisfies Mason's conjecture (I-4).

In light of Example 3.3, for which $\omega$ is it possible to find
a symmetric matrix $A$ with principal
minors $\{\omega(S):S\subseteq E\}$ and off--diagonal entries of $A^{-1}$
equal to $\Cov(X_e,X_f)$  for all $e\neq f$ in $E$?\\

Here are some formulae to keep in mind regarding the covariances for
an arbitrary partition function $Z(\omega;\y)$.
Let $Z^{e}:=Z|_{y_{e}=0}$ and $Z_{e}:=\partial Z/\partial y_{e}$, so that
$Z=Z^{e}+\int \mathrm{d}y_{e}Z_{e}$.  Since $Z$ is multiaffine
$Z=Z^{e}+y_{e}Z_{e}$.  For $e\neq f$
$$\Cov(X_{e},X_{f})= -\frac{y_e y_f}{Z^2}\,\Delta Z\{e,f\}$$
in which
$$\Delta Z\{e,f\}:= Z_{e}Z_{f}- Z_{ef}Z = Z_{e}^{f}Z_{f}^{e}-Z_{ef}Z^{ef}.$$
Also, in general for $e\neq f$,
$$\Cov(X_{e},X_{f})= \frac{\partial^{2}}{\partial y_{e}\partial y_{f}}
\log Z.$$

That convolution of sequences preserves logarithmic concavity (a-2,0)
was shown independently by Karlin \cite{Ka} (on page 394),
by Menon \cite{Men}, and by Hoggar \cite{Hog}.  Liggett \cite{Li}
gives the proof of Lemma 3.5, and moreover proves that the
stronger condition (a-4,0) is preserved.
\begin{LMA}
Let $a_{0}$, $a_{1}$,\ldots and $b_{0}$, $b_{1}$,\ldots be pairwise
commuting indeterminates, and let $b_{-1}=0$.  For each $n\in\NN$ let
$c_{n}:=\sum_{k=0}^{\infty}a_{n+k}b_{k}$.  Then for each $n\geq 1$:
\showon
& & c_{n}^{2}-c_{n-1}c_{n+1}\\
&=& \sum_{k=0}^{\infty}\sum_{j=0}^{k}
\left[ a_{n+j}a_{n+k}-a_{n+j-1}a_{n+k+1} \right]
\left[ b_{j}b_{k}-b_{j-1}b_{k+1} \right].
\showoff
Therefore, for nonnegative real sequences $(a_j)$ and $(b_k)$ such that
every $c_n$ converges, \textup{(a-2,0)} and \textup{(b-2,0)} imply \textup{(c-2,0)}.
\end{LMA}
\begin{proof}
To begin with,
\showon
c_n^2 - c_{n-1} c_{n+1}
&=&
\sum_{j=0}^{\infty}\sum_{k=0}^{\infty}
\left[ a_{n+j} b_j a_{n+k} b_k - a_{n-1+j} b_j a_{n+1+k} b_k \right]\\
&=&
\sum_{0\leq j\leq k}
\left[ a_{n+j} a_{n+k} - a_{n+j-1} a_{n+k+1} \right] b_j b_k\\
& & + \sum_{0\leq k<j}
\left[ a_{n+j} a_{n+k} - a_{n+j-1} a_{n+k+1} \right] b_j b_k.
\showoff
Reindexing the last summation by $k=h-1$ and $j=i+1$, the bounds of
summation are $1\leq h$, $0\leq i$, and $h\leq i+1$, and the general term
of the summand is
$$\left[ a_{n+h-1} a_{n+i+1} - a_{n+h} a_{n+i} \right] b_{h-1} b_{i+1}.$$
When $h=0$ we have $b_{h-1}=b_{-1}=0$, so these terms can be included in the summation.
When $h=i+1$ we have
$$a_{n+h-1} a_{n+i+1} - a_{n+h} a_{n+i} = a_{n+i} a_{n+i+1} - a_{n+i+1} a_{n+i}=0,$$
so these terms can be removed from the summation.  Thus,
\showon
c_n^2 - c_{n-1} c_{n+1}
&=&
\sum_{0\leq j\leq k}
\left[ a_{n+j} a_{n+k} - a_{n+j-1} a_{n+k+1} \right] b_j b_k\\
& & + \sum_{0\leq h\leq i}
\left[ a_{n+h-1} a_{n+i+1} - a_{n+h} a_{n+i} \right] b_{h-1} b_{i+1}\\
&=&
\sum_{0\leq j\leq k}
\left[ a_{n+j} a_{n+k} - a_{n+j-1} a_{n+k+1} \right]
\left[ b_j b_k - b_{j-1} b_{k+1} \right].
\showoff
Under the hypotheses (a-2,0) and (b-2,0) both factors of the general term of
the summand are nonnegative, and (c-2,0) follows.
\end{proof}
\noindent
This extends to doubly--infinite sequences of indeterminates
$a_j$ and $b_k$.  A variety of finiteness or convergence conditions can
then be applied. 

The next result implies the exchangeable case of Conjecture 3.4,
and adds to the equivalent conditions of Theorem 2.7 of Pemantle \cite{Pe}.
\begin{PROP}
For $0\leq k\leq m$, let $e_{k}(\y)$ be the $k$--th elementary 
symmetric function of $\y=\{y_{1},\ldots,y_{m}\}$.  Consider
$$Z(\y)=a_{0}e_{0}(\y)+a_{1}e_{1}(\y)+\cdots+a_{m}e_{m}(\y),$$
in which the $a_{k}$ are nonnegative real numbers.
The following are equivalent:\\
\textup{(a)}\ 
the polynomial $Z(\y)$ satisfies the Rayleigh condition;\\
\textup{(b)}\ 
the sequence $(a_{k}:\ 0\leq k\leq m)$ is logarithmically concave with
no internal zeros, the condition \textup{(a-2,0)};\\
\textup{(c)}\ $\Delta Z\{1,2\}$ is a positive linear combination of
Schur functions;\\
\textup{(d)}\ $\Delta Z\{1,2\}$ is a positive linear combination of
monomial symmetric functions.
\end{PROP}
\begin{proof}
First, to prove that (a) implies (b), assume that $Z(\y)$ satisfies the
Rayleigh condition.  Fix
$1\leq k\leq m-1$.  For $t>0$ let $y_i=t^{-1}$ for $1\leq i\leq k-1$
and $y_j=t$ for $k+2\leq j\leq m$.  Then, as $t\goesto 0$,
$$t^{k-1}Z(\y)\goesto a_{k-1}+a_k(y_{k}+y_{k+1})+a_{k+1}y_{k}y_{k+1}.$$
Since $Z$ satisfies the Rayleigh condition, the inequality
$\Delta Z\{k,k+1\}\geq 0$ for all $t>0$ implies that
$a_k^2\geq a_{k-1}a_{k+1}$.  Proposition 4.7(a) implies that (a-0) holds. 

Next, to prove that (b) implies (c), assume that the sequence $(a_{k}:\ 0\leq k\leq m)$
is logarithmically concave with no internal zeros.
With the notation of Lemma 3.5, if
$b_k=e_k(y_3,...,y_m)$
then
$$Z_1^2(\y)=Z_2^1(\y)=c_1\ \ \mathrm{and}\ \ Z_{12}(\y)=c_2
\ \mathrm{and}\ \ Z^{12}(\y)=c_0.$$
Thus, by Lemma 3.5,
\showon
\Delta Z\{1,2\}
&=& Z_1^2 Z_2^1 - Z_{12} Z^{12} = c_1^2 - c_0 c_2\\
&=&  \sum_{0\leq j\leq k}
\left[ a_{j+1}a_{k+1}-a_{j}a_{k+2} \right]
\left[ b_{j}b_{k}-b_{j-1}b_{k+1} \right].
\showoff
The factor $a_{j+1}a_{k+1}-a_{j}a_{k+2}$ is nonnegative by the
hypothesis of logarithmic concavity with no internal zeros.
The factor $b_{j}b_{k}-b_{j-1}b_{k+1}$ is, by the Jacobi--Trudy
formula, the Schur function of shape $[2^j\ 1^{k-j}]$ in the quantities
$\{y_3,...,y_m\}$.

To see that (c) implies (d), one need only note that Schur functions are positive
linear combinations of monomial symmetric functions.

Finally, if (d) holds then $\Delta Z\{1,2\}$ is nonnegative for all
$\y>\zero$.  Since $Z$ is a symmetric function this suffices to verify (a) that
$Z$ is Rayleigh.
\end{proof}

In general, the difference $\Delta Z\{e,f\}$ might have negative coefficients,
even though $Z$ is Rayleigh.  This happens for example when $Z=Z(\BB\mathsf{K}_4;\y)$
and $e,f\in E(\mathsf{K}_4)$ are non--adjacent edges. See Example 5.2 and Conjectures
5.3 and 5.4.\\

3.3.\ \textsc{supplementary remarks on the big conjecture.}\\

Conjecture 3.4 is implied by the conjunction of two others.
Partition $E=E_{1}\cup E_{2}$ with $E_{1}\cap E_{2}=\none$.
For $i=1,2$ let $\A'_{i}$ be an upward--closed 
subset of $\B(E_{i})$, and let
$$\A_{i}:=\{S\subseteq E:\ S\cap E_{i}\in\A'_{i}\}.$$
This defines two increasing events $\A_{1}$ and $\A_{2}$ with disjoint support.

\begin{CONJ}
Let $Z(\omega;\y)$ satisfy the Rayleigh condition.  Then for all 
$\y>\zero$, and for all pairs of increasing events $\A_{1}$, $\A_{2}$ with 
disjoint support,
$$\la\one_{\A_1\cap \A_2}\ra \leq
\la\one_{\A_1}\ra \cdot \la\one_{\A_2}\ra.$$
\end{CONJ}
This says that pairwise negative correlations (in the strong sense of 
the Rayleigh condition) imply negative correlations for all pairs of
increasing events with disjoint support.  This stronger negative correlation
property is known as \emph{negative association} of the variables $\{X_e:\ 
e\in E\}$, or of their partition function.

The following important result is due to Feder and Mihail \cite{FM};\ also
see Theorem 6.5 of Lyons \cite{Lyo} and Theorem 1.3 of Pemantle \cite{Pe}.

\begin{THM}
If $Z(\omega;\y)$ is Rayleigh and homogeneous then it is
negatively associated for all $\y>\zero$.
\end{THM}

\begin{CONJ}
Let $Z(\omega;\y)$ be negatively associated for all $\y>\zero$.
Then $(f_{k}(\omega)/\binom{m}{k}:\ 0\leq k\leq m)$ is 
logarithmically concave with no internal zeros, the condition
\textup{(f-4,0)}.
\end{CONJ}

Allowing the $\y>\zero$ to vary independently is essential
for Conjecture 3.4, as the following example shows.
\begin{EG}\emph{
For $\gamma\geq 0$, let 
$$Z(t) = 1 + 12t + 60t^2 + 20\gamma t^3 + 60t^4 + 12t^5 + t^6.$$
This sequence of coefficients $(f_k)$ satisfies (f-4) if and only if
$4\leq\gamma\leq 8$, (f-2) if and only if $3\leq\gamma\leq 15 $,
(f-1) if and only if $\gamma\geq 3$, and (f-0) if and only if $\gamma>0$.
For an exchangeable set of binary variables $\{X_1,\ldots,X_{6}\}$
with this partition function $Z(t)$, all pairwise correlations are the
same as for $X_{1}$ and $X_{2}$.  A short calculation yields that
if $2.61725\approx\sqrt{6.85} \leq\gamma\leq 8$ then $\Delta Z\{1,2\}$
is a polynomial with nonnegative
coefficients, and hence is nonnegative for all $t>0$.  Thus, for $\gamma$ in this
range $Z(t)$ satisfies a restricted form of the Rayleigh condition (the form in which
all $y_e=t$ are equal).  But when $\gamma=2.9$, for example, the unimodality
condition (f-1) does not hold.  Thus, in Conjecture 3.4, it is essential that
the $\y>\zero$ be allowed to vary at least somewhat independently.
}\end{EG}

\vspace{4mm}

3.4.\ \textsc{an equivalent form of the big conjecture.}\\

To reduce the general case to the exchangeable case consider the
\emph{symmetrizing operator} $Z\mapsto \zZ$, acting on a partition function
$Z(\omega;\y)$ of $m$ variables $y_1,...,y_m$ by
$$  Z(\widetilde{\omega};\y):=\zZ(\omega;\y):=
\frac{1}{m!}\sum_{\sigma\in\S_m}Z(\omega;\y_\sigma),$$
the sum being over all permutations $\sigma$ in the symmetric group $\S_m$
(on the set $E$ indexing $\y$), and
$\y_\sigma^S:=\prod_{e\in S} y_{\sigma(e)}$
for all $S\subseteq E$.
Calculation reveals that for all $S\subseteq E$
with $|S|=k$, $\widetilde{\omega}(S)=f_k(\omega)/\binom{m}{k}$.  Thus
$$\zZ(\omega;\y) = \sum_{k=0}^m \frac{f_k(\omega)}{\binom{m}{k}} e_k(\y).$$
Notice that $f_k(\widetilde{\omega})=f_k(\omega)$ for all $0\leq k\leq m$.

\begin{CONJ}
If $Z(\omega;\y)$ is Rayleigh then $\zZ(\omega;\y)$ is Rayleigh.
\end{CONJ}
Proposition 3.6 and elementary logic establish the following.
\begin{PROP}
Conjecture $3.11$ is equivalent to Conjecture $3.4$.
\end{PROP}
We are faced with the challenge of proving Conjecture 3.11, so far without
success.  Here are some concrete ideas towards a proof, and a reasonably
plausible sufficient condition.

Assume that $Z=Z(\omega;\y)$ is Rayleigh.  To show that $\zZ$ is
Rayleigh it suffices to show that $\Delta\zZ\{1,2\}\geq 0$ for all
$\y>\zero$, since $\zZ$ is exchangeable (a symmetric function).
By Proposition 3.6, this happens if and only if $\Delta\zZ\{1,2\}$
has positive coefficients as a polynomial in $\y$.

For a symmetric function $F$ with $m-k$
indeterminates and $S\subseteq \{1,...,m\}$ with $|S|=k$, let
$\natural_S F$ denote the same symmetric function
of the indeterminates $\{y_i:\ i\in \{1,...,m\} \drop S\}$.
$$
\zZ_1 = \frac{1}{m}\sum_{e\in E}\natural_1(Z_e)^\sim.
$$
$$
\zZ_{12} =
\frac{1}{\binom{m}{2}}\sum_{\{e,f\}\subseteq E}  \natural_{12} (Z_{ef})^\sim.
$$
For any symmetric function $F$ of $m-1$ indeterminates,
$$(\natural_1 F)|_{y_2=0} = (\natural_2 F)|_{y_1=0}.$$
For each $e\in E$ let
$$ z(e):= (\natural_1 (Z_e)^\sim)^2 = (\natural_2 (Z_e)^\sim)^1.$$
\showon
\zZ_1^2 \zZ_2^1
&=&
\left(\frac{1}{m}\sum_{e\in E}z(e)\right)
\left(\frac{1}{m}\sum_{e\in E}z(e)\right)\\
&=& \frac{1}{m^2} \sum_{\{e,f\}\subseteq E}
\left[2z(e)z(f)+\frac{1}{m-1}(z(e)^2+z(f)^2)\right]\\
&\geq& \frac{1}{m^2}
\frac{2m}{m-1}\sum_{\{e,f\}\subseteq E}z(e)z(f).\\
\zZ_{12}\zZ^{12}
&=&
\frac{2}{m(m-1)}\sum_{\{e,f\}\subseteq E}  (\natural_{12} (Z_{ef})^\sim) \zZ^{12}.\\
\Delta\zZ\{1,2\}
&\geq& \frac{1}{\binom{m}{2}}
\sum_{\{e,f\}\subseteq E}(z(e)z(f)-z(ef)\zZ^{12}),
\showoff
in which $z(ef):=\natural_{12} (Z_{ef})^\sim$.  Thus $\Delta\zZ\{1,2\}$ is bounded
below by an average of $\binom{m}{2}$ terms, one for each $\{e,f\}\subseteq E$.

\begin{CONJ}
If $Z(\omega;\y)$ is Rayleigh then 
$$\sum_{\{e,f\}\subseteq E} (z(e)z(f) - z(ef)\zZ^{12})\geq 0$$
for all $\y>\zero$.
\end{CONJ}
By the preceding calculations, Conjecture 3.13 implies Conjecture 3.11. 
Here is a simpler and stronger conjecture.
\begin{CONJ}
If $Z(\omega;\y)$ is Rayleigh then for all $\{e,f\}\subseteq E$,
$$z(e)z(f) - z(ef)\zZ^{12}\geq 0$$
for all $\y>\zero$.
\end{CONJ}

\section{Necessary Conditions.}

The hypothesis that $Z(\omega;\y)$ satisfies the Rayleigh condition imposes some
strong -- and perhaps surprising -- necessary conditions on the collection of sets
$S\subseteq E$ for which $\omega(S)>0$.  It also implies some inequalities on
the values of $\omega$ that are an all--pervasive local form of logarithmic concavity.

Given $\omega:\B(E)\goesto[0,\infty)$, we let
$$\Supp(\omega):=\{S\subseteq E:\ \omega(S)>0\},$$
and call this the \emph{support} of $\omega$, or of $Z(\omega;\y)$.
This section is mainly concerned with the combinatorial properties of
$\Supp(\omega)$ when $Z(\omega;\y)$ is Rayleigh.\\

4.1.\ \textsc{convexity and logarithmic submodularity.}\\

If $\Q$ is the support of $Z$ then the support of $Z^e$ is
$$\Q^e:=\{S:\ e\not\in S\ \mathrm{and}\ S \in \Q \},$$
and the support of $Z_e$ is
$$\Q_e:=\{S\drop\{e\}:\ e\in S\ \mathrm{and}\ S \in \Q\}.$$
This notation is extended to multiple (distinct) indices in the obvious way.
The \emph{dual} of a set--system $\Q\subseteq \B(E)$ is
$$\Q^*:=\{E\drop S:\ S\in\Q\}.$$

A set--system $\Q$ is \emph{Rayleigh} if the polynomial
$Z(\Q;\y)$ of Example 3.1 meets the Rayleigh condition.
A set--system $\Q$ is \emph{weakly Rayleigh} if there is some nonnegative weight
function $\omega$ with support equal to $\Q$ such that
$Z(\omega;\y)$ meets the Rayleigh condition. 

\begin{LMA}
Let $\Q\subseteq\B(E)$ be a (weakly) Rayleigh set--system.\\
\textup{(a)}\ For all $g\in E$, both $\Q^g$ and $\Q_g$ are (weakly) Rayleigh.\\
\textup{(b)}\ The dual $\Q^*$ is (weakly) Rayleigh.
\end{LMA}
\begin{proof} 
Let $Z(\omega;\y)$ be a polynomial with support $\Q$.  To prove part (a),
for distinct $e,f,g\in E$, a short calculation shows that
$$\Delta Z\{e,f\}=
\Delta Z^g\{e,f\}
+y_g\Theta Z\{e,f|g\}
+y_g^2 \Delta Z_g\{e,f\}$$
for some polynomial $\Theta Z\{e,f|g\}$.  Taking limits as $y_g\goesto 0$
or as $y_g\goesto\infty$ shows that if $Z$ meets the Rayleigh condition
then so do $Z^g$ and $Z_g$.  If $\omega=\one_\Q$ then $\Q^g$ and $\Q_g$ are
Rayleigh as well.

For part (b) one calculates that for $e,f\in E$,
$$\Delta Z^*\{e,f\}(\y)= (\y^E)^2\Delta Z\{e,f\}(1/\y),$$
from which the result follows.
\end{proof}

\begin{THM}
Let $\Q\subseteq \B(E)$ be a weakly Rayleigh set--system.
If $\{\none,E\}\subseteq\Q$ then $\Q=\B(E)$.
\end{THM}
\begin{proof}
Let $Z(\omega;\y)$ be a polynomial with support $\Q$ that meets the
Rayleigh condition.

We proceed by induction on $|E|=m$.  The bases $m=0$
or $m=1$ are trivial.  For the case $m=2$ let $E=\{e,f\}$, so that $Z$
has the form
$$Z= A+ B y_e + C y_f + D y_e y_f$$
for nonnegative constants $A,B,C,D$ with $A$ and $D$ positive.
Since $Z$ is Rayleigh the inequality $BC\geq AD$ holds,
so that both $B$ and $C$ are positive as well.  Hence $\Q=\B(E)$ in this case.

For the induction step we assume that $m\geq 2$.  Arguing for a contradiction,
suppose that $S\subset E$ is such that $S\not\in\Q$.  (So, in particular,
$S\not\in\{\none,E\}$.)  If there is a subset
$\none\subset T\subset S$ such that $T\in\Q$ then the interval $[T,E]$ of
$\B(E)$ is a smaller Boolean algebra, and the contraction $\Q_T$ is weakly
Rayleigh by Lemma 4.1(a).  By the induction hypothesis, $S\drop T\in \Q_T$,
so that $S\in\Q$, a contradiction.  Therefore, $\Q\cap[\none,S]=\{\none\}$.
Thus, there is an element $e\in S$ such that $\{e\}\not\in\Q$.

Now, if $e\in T\subset E$ then $T\not\in\Q$;  for if it were the case that
$T\in\Q$ then the fact that the deletion $\Q^{E\drop T}$ is Rayleigh by Lemma 4.1(a)
and the induction hypothesis imply that $\{e\}\in\Q$.  Thus, $\Q\cap[\{e\},E]=\{E\}$.
Therefore $Z(\omega;\y)$ has the form $Z=Z^e+K\y^E$ for some constant $K>0$.
Since $m\geq 2$ there is some $f\in E\drop\{e\}$.  Now $\Q_e^f=\none$ and
both $\Q_{ef}\neq\none$ and $\Q^{ef}\neq\none$, so that $\Delta Z\{e,f\}<0$
for \emph{every} $\y>\zero$.  This contradicts the hypothesis that $Z(\omega;\y)$
is Rayleigh, completing the induction step and the proof.
\end{proof}

A set--system $\Q$ is \emph{convex} if it satisfies the condition that
for any $S\subseteq T\subseteq S'\subseteq E$, if $S,S'\in\Q$ then $T\in\Q$.

\begin{CORO}
Every weakly Rayleigh set--system  is convex.
\end{CORO}
\begin{proof} 
Let $\Q\subseteq \B(E)$ be a weakly Rayleigh set--system.
If $S\subseteq T\subseteq S'$ with $S,S'\in\Q$ then consider the set--system
$(\Q_S)^{E\drop S'}$.  This is weakly Rayleigh by Lemma 4.1(a), and contains both
$\none$ and $S'\drop S$.  By Theorem 4.2, $T\drop S$ is in $(\Q_S)^{E\drop S'}$,
so that $T$ is in $\Q$.
\end{proof}

\begin{THM}
If $\omega:\B(E)\goesto[0,\infty)$ is Rayleigh
then $\omega$ is \emph{logarithmically submodular:}\
for all $S,T\in\B(E)$,
$$\omega(S)\omega(T)\geq\omega(S\cap T)\omega(S\cup T).$$
\end{THM}
\begin{proof} 
By the $m=2$ case of the proof of Theorem 4.2, this inequality
holds whenever $S\cap T$ is covered by both $S$ and $T$.
It holds trivially if either of
$S\cap T$ or $S\cup T$ is not in $\Supp(\omega)$, so assume otherwise.
By Corollary 4.3, the interval $[S\cap T,S\cup T]$ is contained in $\Supp(\omega)$.
Let
$$S\cap T=S_0\subset S_1\subset\cdots\subset S_k=S$$ and
$$S\cap T=T_0\subset T_1\subset\cdots\subset T_\ell=T$$
be saturated chains in $\B(E)$.  By the above remarks, for all $1\leq i\leq k$
and $1\leq j\leq\ell$,
$$\omega(S_i\cup T_{j-1})\omega(S_{i-1}\cup T_j)
\geq
\omega(S_{i-1}\cup T_{j-1})\omega(S_i\cup T_i).$$
Taking the product of all these inequalities and cancelling the common
factors (which are strictly positive),
we obtain
$$\omega(S_k\cup T_0)\omega(S_0\cup T_\ell)\geq
\omega(S_0\cup T_0)\omega(S_k\cup T_\ell).$$
That is
$$\omega(S)\omega(T)\geq\omega(S\cap T)\omega(S\cup T),$$
as desired.
\end{proof}
Interestingly, the choice of saturated chains in this proof is arbitrary
and disappears in the answer.

\vspace{4mm}

4.2.\ \textsc{exchange properties of the support.}\\

\begin{LMA}
Let $\Q\subseteq\B(E)$ be a weakly Rayleigh set--system, and let $A,B\in\Q$
with $A\cap B=\none$.  For every $\{e,f\}\subseteq B$ and $g\in A$,
at least one of the sets $\Q_{eg}^f$ or $\Q_{fg}^e$ is not empty.
\end{LMA}
\begin{proof} 
Let $Z(\omega;\y)$ be a Rayleigh polynomial with support $\Q$,
and suppose that the conclusion fails to hold.
Let $\{e,f\}\subseteq B$ and $g\in A$ be such that both
$\Q_{eg}^f$ and $\Q_{fg}^e$  are empty.  Then, in $\Delta Z\{e,f\}$
the indeterminate $y_g$ does not occur in the term $Z_e^f Z_f^e$.
However, the pair $(A,B)$ contributes $\y^A \y^B$ to the term
$Z^{ef}Z_{ef}$, so the indeterminate $y_g$ occurs in this term.
No matter what values $y_c>0$ are chosen for all $c\neq g$, as
$y_g\goesto\infty$, $\Delta Z\{e,f\}\goesto -\infty$.  This contradicts
the hypothesis that $Z$ meets the Rayleigh condition, completing the proof.
\end{proof}

A \emph{delta--matroid} is a set--system $\Q\subseteq\B(E)$ that 
satisfies the following \emph{symmetric exchange axiom}:\\
\textbf{(SEA)}\  if $A,B\in\Q$ and $e\in A\tri B$, then there is an
$f\in A\tri B$ such that $A\tri\{e,f\}\in\Q$.\\
(Here $\tri$ denotes the symmetric difference of sets.)  Notice that $\Q$
is a delta--matroid if and only if $\Q^*$ is a delta--matroid.
A good deal of matroid theory generalizes well to delta--matroids;\ see
\cite{B1,BC,BJ,FT,GIM} for starters.

\begin{THM}
Every weakly Rayleigh set--system is a convex delta--matroid.
\end{THM}
\begin{proof}
Let $Z(\omega;\y)$ be a polynomial with support $\Q$ that meets the
Rayleigh condition.  That $\Q$ is convex has been established in Corollary 4.3.
The strategy of the proof is along the lines of \cite{Choe1,COSW}.

We verify that $\Q$ is a delta--matroid by induction on the
size $|E|=m$ of the underlying set.  The base cases $m\leq 2$ are trivial,
so assume that $m\geq 3$.  By Lemma 4.1(a) and the induction hypothesis, 
for any $g\in E$ the set systems $\Q_{g}$ and $\Q^{g}$ are convex 
delta--matroids.

Now consider any $A,B\in\Q$ and $e\in A\tri B$.
To verify the SEA for $(A,B,e)$ in $\Q$, we must find an element
$f\in A\tri B$ such that $A\tri \{e,f\}\in\Q$.
If $A\tri\{e\}\in\Q$ then we can choose $f=e$ to satisfy the SEA,
so we are left with the case that $A\tri\{e\}\not\in\Q$.

If $g\in A\cap B$ then consider the sets $A':=A\drop\{g\}$ and $
B':=B\drop\{g\}$ in $\Q_g$, and the element $e\in A'\tri B'$.
By the SEA for $(A', B', e)$ in $\Q_{g}$, there is an element $f\in A'\tri B'$
such that $A'\tri\{e,f\}\in\Q_g$.  This is an element $f\in A\tri B$
such that $A\tri\{e,f\}\in\Q$.  Thus, we can assume that $A\cap B=\none$.

If $g\in E\drop(A\cup B)$ then consider the sets $A$ and $B$ in $\Q^g$,
and the element $e\in A\tri B$.   By the SEA for $(A,B,e)$ in $\Q^{g}$,
there is an element $f\in A\tri B$ such that $A\tri\{e,f\}\in\Q^g$.  This
is an element $f\in A\tri B$ such that $A\tri\{e,f\}\in\Q$.  Thus, we can
assume that $A\cup B=E$.

Now, if $A=\none$ then $B=E$, while if $B=\none$ then $A=E$.  In either case,
since $\Q$ is convex, $\Q=\B(E)$.  Since $\B(E)$ is a delta--matroid we can
assume that both $A$ and $B$ are nonempty.  Since $|E|=m\geq 3$, one of the sets
$A$ or $B$ has at least two elements.  By Lemma 4.5 there is a third set $C\in\Q$
such that $A\cap C\neq\none$ and $B\cap C\neq\none$.

From this point on we resort to a case analysis based on the two main cases
$e\in A$ or $e\in B$, and on several subcases.

\noindent\textbf{Case I:}\ $e\in A$.

\underline{Subcase (i):}\ $e\in A\drop C$.\\
Let $g\in A\cap C$, and consider $A':=A\drop\{g\}$
and $C':=C\drop\{g\}$ in $\Q_g$.  By the SEA for $(A',C',e)$ in 
$\Q_{g}$, there is an $f\in A'\tri C'$ such that $A'\tri\{e,f\}\in\Q_g$.
Since $A'\tri C' \subseteq A\tri B$, this is an element $f\in A\tri B$ such
that $A\tri\{e,f\}\in\Q$.  This verifies the SEA for $(A,B,e)$ in $\Q$
in this subcase.

\underline{Subcase (ii):}\ $e\in A\cap C$ and $|A\cap C|\geq 3$.\\
Let $g\in B\cap C$, so that $C':=C\drop\{g\}$ and $B':=B\drop\{g\}$ are in $\Q_g$.
By the SEA for $(C',B',e)$ in $\Q_{g}$, there
is an $f\in C'\tri B'$ such that $C'\tri\{e,f\}$ is in $\Q_g$.  Thus,
$C'':=C\tri\{e,f\}$ is in $\Q$.  Now $e\in A\drop C''$, and $g\in B\cap C''$,
and since $|A\cap C|\geq 3$ it follows that $|A\cap C''|\geq 1$.  Thus, $C''$ is a
set with the properties of $C$ in subcase I(i), reducing subcase I(ii) to
that previously solved subcase.

\underline{Subcase (iii):}\ $e\in A\cap C$ and $|A\cap C|=2$.\\
Repeating the argument for subcase I(ii) produces a set $C''$ with 
the properties of $C$ in subcase I(i) except when $A\cap C=\{e,f\}$, in which case
$C'':=C\tri\{e,f\}$ is disjoint from $A$.  But then
$[C'',C]\subseteq \Q$ since $\Q$ is convex, and
it follows that $C''':=C\drop\{e\}\in\Q$.
Now this $C'''$ is a set like $C$ in subcase I(i), reducing subcase I(iii) to
that previously solved subcase.

\underline{Subcase (iv):}\ $A\cap C=\{e\}$.\\
If $|A|\geq 2$ and $|B|\geq 2$ then this case can be avoided as follows:\
for any $g\in A\drop\{e\}$,  Lemma 4.5 can be used to ensure that the
particular element $g$ is contained in $C$, so that one of the subcases
I(i,ii,iii) holds instead.  Thus, we can assume that either $|A|=1$ or
$|B|=1$.

If $A=\{e\}$ then $|B|\geq 2$, and $[\{e\},C]\subseteq\Q$ since $\Q$ is convex.
Thus we can assume that $|C|=2$, so let $C=\{e,g\}$ with $g\in B$.  If
$\{g\}\in\Q$ then $A\tri\{e,g\}=\{g\}\in\Q$ suffices to verify the
SEA for $(A,B,e)$ in $\Q$, so we can assume that $\{g\}\not\in\Q$.
Applying the SEA to $(\{e\},B\drop\{g\},e)$ in $\Q_g$, there is an
$f\in B\drop\{g\}$ such that $\{f,g\}\in\Q$.   If $\{f\}\in\Q$ then
$A\tri\{e,f\}=\{f\}\in\Q$ suffices to verify the
SEA for $(A,B,e)$ in $\Q$, so we can assume that $\{f\}\not\in\Q$.
Now let $\R:=\Q\cap\B(\{e,f,g\})$.  The partition polynomial for the
restriction of $\omega$ to $\R$ has the form
$$Z|_\R=K_e y_e + K_{ef} y_e y_f + K_{eg} y_e y_g + K_{fg} y_f y_g + K_{efg}
y_e y_f y_g$$
for nonnegative constants $K_{S}$ with $K_e$, $K_{eg}$, and $K_{fg}$
positive.  (Note that $\none\not\in\R$ since $\none=A\tri\{e\}\not\in\Q$.)
This is Rayleigh, by Lemma 4.1(a).
The Rayleigh difference of $f$ and $g$ in $Z|_\R$ is
\showon
\Delta Z|_\R\{f,g\}
&=& K_{ef} K_{eg} y_e^2 - K_e y_e (K_{fg}+ K_{efg} y_e)\\
&=& (K_{ef} K_{eg}-K_e K_{efg})y_e^2 - K_e K_{fg} y_e.
\showoff
This quantity can be made negative by taking $y_e>0$ to be sufficiently
small, contradicting the fact that $Z|_\R$ is Rayleigh.
Thus this part of the subcase does not arise, completing the proof of
subcase I(iv) when $|A|=1$.

If $|B|=1$ then $|A|\geq 2$. Let $B=\{g\}$, so that $C=\{e,g\}$.
Applying the SEA to $(\{g\},A\drop\{e\},g)$ in $\Q_e$,
either $\{e\}\in\Q$ or there is an $h\in A\drop \{e\}$ such that
$\{e,h\}\in\Q$.  If $\{e\}\in\Q$ then $[\{e\},A]\subseteq\Q$ since $\Q$ is convex;\
thus we can assume that there is an $h\in A\drop \{e\}$ such that $\{e,h\}\in\Q$.
If $\{h\}\in\Q$ then, since $\Q$ is convex, $A\tri\{e\}\in\Q$;\ since
we have assumed that this is not the case, $\{h\}\not\in\Q$.
Now let $\R:=\Q\cap\B(\{e,g,h\})$.  The partition polynomial for the
restriction of $\omega$ to $\R$ has the form
$$Z|_\R=K_e y_e + K_g y_g + K_{eg} y_e y_g + K_{eh} y_e y_h + K_{gh} y_g y_h
 + K_{egh} y_e y_g y_h$$
for nonnegative constants $K_{S}$ with $K_e$, $K_g$, $K_{eg}$, and $K_{eh}$
positive.
This is Rayleigh, by Lemma 4.1(a).
The Rayleigh difference of $e$ and $h$ in $Z|_\R$ is
$$\Delta Z|_\R\{e,h\}=
(K_{e}+K_{eg}y_g)K_{gh}y_g - K_{g} y_{g}(K_{eh}+ K_{egh} y_g).$$
Since this is nonnegative for all $\y>\zero$ it must be the case that $K_{gh}>0$;\
that is, $\{g,h\}\in\Q$.  This set has the form of the set $C$ in subcase I(i),
reducing the problem to that previously solved subcase.  This completes the analysis
of Case I.

\noindent\textbf{Case II:}\ $e\in B$.

\underline{Subcase (i):}\ $|A|\geq 2$.\\
In this case Lemma 4.5 can be used to ensure that the particular element
$e\in B$ is contained in $C$.  From this point, the argument for subcase
I(i) applies verbatim to this subcase as well, establishing the
SEA for $(A,B,e)$ in $\Q$ in this subcase.

\underline{Subcase (ii):}\ $|A|=1$.\\
If $e\in C$ then $A\subset A\tri\{e\}\subseteq C$, so that $A\tri\{e\}\in\Q$
because $\Q$ is convex.  Thus we can assume that $e\in B\drop C$.
Let $A=\{g\}$.  Since $\{g\}\subseteq C$ and $\Q$ is convex we can 
assume that $|C|=2$, so let $C=\{g,h\}$ with $h\in B$.  By the SEA
for $(\{g\},B\drop\{h\},e)$ in $\Q_{h}$, either $\{e,g,h\}\in\Q$ or
$\{e,h\}\in\Q$.  If $\{e,g,h\}\in \Q$ then $\{e,g\}\in\Q$ since $\Q$ 
is convex;\ but $\{e,g\}=A\tri\{e\}\not\in\Q$ by a previous assumption.
Thus we see that $\{e,h\}\in\Q$.  Let $\R:=\Q\cap\B(\{e,g,h\})$.
The partition polynomial for the restriction of $\omega$ to $\R$ has the form
$$Z|_\R = K_{\none} + K_e y_e + K_g y_g + K_{h} y_{h}
+ K_{eh} y_e y_h + K_{gh} y_g y_h$$
for nonnegative constants $K_{S}$ with $K_g$, $K_{gh}$, and $K_{eh}$
positive.  This is Rayleigh, by Lemma 4.1(a).
The Rayleigh difference of $e$ and $h$ in $Z|_\R$ is
$$\Delta Z|_{\R}\{e,h\}=
K_{e}(K_{h}+K_{gh}y_g) - (K_{\none}+K_{g} y_{g})K_{eh}.$$
Since this is nonnegative for all $\y>\zero$ it must be the case that
$K_{e}>0$; that is, $\{e\}\in\Q$.  Since $A\tri\{e,g\}=\{e\}\in\Q$,
this establishes the SEA for $(A,B,e)$ in $\Q$ in Case II.

This completes the case analysis, the induction step, and the proof.
\end{proof}

After that, proofs of the following consequences follow familiar lines.
Proposition 4.7(a) implies (f-0) for Rayleigh $\omega$ with support $\Q$. 

\begin{PROP}
Let $\Q$ be a convex delta--matroid, and let $A,B\in\Q$.\\ 
\textup{(a)}\  If $|A|<|B|$ then there exists $b\in B\drop A$ such that
$A\cup\{b\}\in\Q$.\\
\textup{(b)}\  If $|A|<|B|$ then there exists $b\in B\drop A$ such that
$B\drop\{b\}\in\Q$.\\
\textup{(c)}\  If $|A|=|B|$ then for every $a\in A\drop B$
there is a $b\in B\drop A$ such that $A\drop\{a\}\cup\{b\}\in\Q$.\\
\textup{(d)}\  If $|A|=|B|$ then for every $a\in A\drop B$
there is a $b\in B\drop A$ such that $B\drop\{b\}\cup\{a\}\in\Q$.
\end{PROP}

\begin{CORO}
If $\Q$ is weakly Rayleigh then all maximal elements of $\Q$ have
the same cardinality $r$, and all minimal elements of $\Q$ have the same
cardinality $s$.
\end{CORO}

\begin{CORO}
Let $\Q$ be a homogeneous set--system.  If $\Q$ is weakly Rayleigh
then $\Q=\BB\M$ is the set of bases of a matroid $\M$.
\end{CORO}

\begin{CORO}
Let $\Q$ be a simplicial complex.  If $\Q$ is weakly Rayleigh
then $\Q=\II\M$ is the set of independent sets of a matroid $\M$.
\end{CORO}
Corollary 4.10 is bad news for the conjectures in Section 2.2.
The Big Conjecture 3.4 is directly relevant only to Mason's Conjecture (I-4).

\begin{CORO}
If $\Q$ is (weakly) Rayleigh then $\Q=\II\M\cap\SS\N$ is the intersection of 
the set of independent sets of a matroid $\M$ and
the set of spanning sets of a matroid $\N$.  Moreover,
both $\BB\M$ and $\BB\N$ are (weakly) Rayleigh.
\end{CORO}

It is natural to wonder:\ does every convex delta--matroid have the LYM property?\\

4.3.\ \textsc{flattening a convex delta--matroid.}\\

Let $\Q\subseteq\B(E)$, let $r:=\max\{|S|:\ S\in\Q\}$ and
$s:=\min\{|S|:\ S\in\Q\}$, and let $\ell:=r-s$.  Assume that $E\cap\{1,...,\ell\}=\none$.
Let $E^\flat:=E\cup\{1,...,\ell\}$, and define $\Q^\flat\subseteq\B(E^\flat)$ by
$$\Q^\flat:=\{S\subseteq E^\flat:\ |S|=r\ \mathrm{and}\ S\cap E\in\Q\}.$$

Given $\omega:\B(E)\goesto[0,\infty)$, define
$\omega^\flat:\B(E^\flat)\goesto[0,\infty)$ by putting
$$\omega^\flat(S):=
\left\{\begin{array}{ll}
\omega(S\cap E) & \mathrm{if}\ |S|=r,\\
0 & \mathrm{if}\ |S|\neq r.
\end{array}\right.$$
for every $S\subseteq E^\flat$.  It follows that the support
of $\omega^\flat$ is $\Supp(\omega)^\flat$.  The partition function is
$$Z(\omega^{\flat};\y)=\sum_{S\subseteq 
E}\omega(S)\y^{S}e_{r-|S|}(y_{1},\ldots,y_{\ell}).$$

\begin{THM}
Let $\Q$ be a convex delta--matroid.  Then $\Q^\flat$ is the set
of bases of a matroid.
\end{THM}
\begin{proof}
To verify the basis exchange axiom for $\Q^\flat$, consider any $A,B\in\Q^\flat$
and $a\in A\drop B$.  Let $A':=A\cap E$, $A'':= A\drop E$, $B':=B\cap E$, and
$B'':= B\drop E$.  There are two main cases:\ either $a\in A'\drop B'$ or
$a\in A''\drop B''$.

If $a\in A'\drop B'$ then either $A'\drop\{a\}\not\in\Q$ or $A'\drop\{a\}\in\Q$.
If $A'\drop\{a\}\not\in\Q$ then, by the SEA applied to $(A',B',a)$ in $\Q$
and since $\Q$ is convex, there is a $b\in B'$ such that $A'\cup\{b\}\drop\{a\}\in\Q$.
Thus, $A\cup\{b\}\drop\{a\}\in\Q^\flat$.  On the other hand, if $A'\drop\{a\}
\in\Q$ then either $|A'|\leq |B'|$ or $|A'|>|B'|$.  If $|A'|\leq |B'|$ then
by Proposition 4.7(a) there is a $b\in B'\drop A'$ such that $A'\cup\{b\}\drop\{a\}
\in\Q$, so that $A\cup\{b\}\drop\{a\}\in\Q^\flat$.  If $|A'|>|B'|$ then
$B''\drop A''\neq\none$, and for any $b\in B''\drop A''$ we have
$A\cup\{b\}\drop\{a\}\in\Q^\flat$.

If $a\in A''\drop B''$ then either $|A''|\leq|B''|$ or $|A''|>|B''|$.
If $|A''|\leq|B''|$ then $a\in A''\drop B''$ implies that $B''\drop A''\neq\none$;
for any $b\in B''\drop A''$ we have $A\cup\{b\}\drop\{a\}\in\Q^\flat$.
If $|A''|>|B''|$ then $|A'|<|B'|$, so by Proposition 4.7(a) there $b\in B'\drop A'$
such that $A'\cup\{b\}\drop\{a\}\in\Q$, so that $A\cup\{b\}\drop\{a\}\in\Q^\flat$.
This verifies the matroid basis exchange axiom for $\Q^\flat$.
\end{proof}

A \emph{strong map} $\M\goesto\M'$ of matroids $\M$ and $\M'$ is a matroid
$\N$ with a distinguished subset $S\subseteq E(\N)$ such that
$\M\simeq \N\drop S$ and $\M'\simeq \N/S$.
\begin{CORO}
Let $\Q$ be a convex delta--matroid.  For each $s\leq k\leq r$ let $\Q_k$
denote the collection of sets in $\Q$ of size $k$.  Then each $\Q_k=\BB\M_k$
is the set of bases of a matroid $\M_k$, and there are strong maps
$$\M_r \goesto \M_{r-1} \goesto \cdots \goesto \M_{s+1} \goesto \M_s,$$
every composition of which is also a strong map.
\end{CORO}
\begin{proof}
Let $\N$ be the matroid with $\BB\N=\Q^\flat$, and let $\ell:=r-s$.
For each $0\leq j\leq \ell$ let
$\M_{s+j}:=(\N\drop\{1,...,j\})/\{j+1,...,\ell\}$.  These matroids are
such that $\Q_k=\BB\M_k$ for all $s\leq k\leq r$.
For $s\leq j< k\leq r$, let $D:=\{1,...,j-s\}$ and $C:=\{k-s+1,...,r-s\}$
and $S:=\{j-s+1,...,k-s\}$.  The matroid $(\N\drop D)/C$ with distinguished
subset $S$ provides a strong map $\M_k\goesto\M_j$.
\end{proof}

\begin{PROP}
For $\omega:\B(E)\goesto[0,\infty)$, the sequence
$(f_{k}(\omega):\ s\leq k\leq r)$ satisfies \textup{(f-2,0)}
if and only if $Z(\omega^\flat;\y)|_{y_e=1:\ e\in E}$ is Rayleigh.
\end{PROP}
\begin{proof}
Let $Z(y_1,...,y_\ell)$ be the polynomial obtained from $Z(\omega^\flat;\y)$
by setting $y_e=1$ for all $e\in E$.  Then
$$Z(y_1,...,y_\ell)=\sum_{j=0}^\ell f_{s+j}(\omega) e_j(y_1,...,y_\ell).$$
Thus it follows from Proposition 3.6.
\end{proof}

4.4.\ \textsc{the triple condition.}\\

For any $Z=Z(\omega;\y)$ and $\{e,f,g\}\subseteq E$,
$$\Delta Z\{e,f\}=
\Delta Z^g\{e,f\}
+y_g \Theta Z\{e,f|g\}
+y_g^2 \Delta Z_g\{e,f\}$$
in which 
$$\Theta Z\{e,f|g\}:=
Z_e^{fg} Z_{fg}^e + Z_f^{eg} Z_{eg}^f
- Z_g^{ef} Z_{ef}^g - Z_{efg} Z^{efg}.$$

Proposition 4.15 is the analogue of Corollary 3.3 of \cite{CW} in a
more general setting.
\begin{PROP}
Assume that $Z=Z(\omega;\y)$ is Rayleigh, and let
$\{e,f,g\}\subseteq E$.  Then for all $\y>\zero$,
$$\Theta Z\{e,f|g\}\geq -2\sqrt{\Delta Z^g\{e,f\}\Delta Z_g\{e,f\}}.$$
\end{PROP}
\begin{proof}
The quantity $\Delta Z\{e,f\}= C + B y_g + A y_g^2$ is a quadratic
polynomial of $y_g$.  For any values $y_c>0$ for all $c\in E\drop\{e,f,g\}$,
this polynomial is nonnegative for all $y_g\geq 0$.  Since $AC\geq 0$,
if $B^2-4AC\geq 0$ then both roots have the same sign.  Thus, either
$B^2-4AC<0$ or $B^2-4AC\geq 0$ and $B\geq 0$.  This yields the desired
inequality.
\end{proof}

\section{Sufficient Conditions.}

5.1. \ \textsc{examples.}\\

\begin{EG}\emph{
Let $\U=\U(m,r)$ denote the uniform matroid of rank $r$ on $m$ elements.
Its rank function on $\B(E)$ is given by $\mathrm{rank}(S)=\min\{r,|S|\}$
for all $S\subseteq E$.  Its Potts model partition function
$Z(\U,q;\y)$ is in the exchangeable case, and so
it is Rayleigh if and only if $(q^{-\min\{r,k\}}:\ 0\leq k\leq m)$
is logarithmically concave with no internal zeros, by Proposition 3.6.
This occurs for all $0<q\leq 1$.
}\end{EG}

\begin{EG}\emph{
Small graphs can be shown to be $\II$--Rayleigh by moderate computations.
For $Z=Z(\II\mathsf{K}_n;\y)$, for example, there are two cases by symmetry
for $\Delta Z\{e,f\}$:\ the edges $\{e,f\}$ are either adjacent or not.  With
$E(\mathsf{K}_4)=\{\{1,2\},...,\{3,4\}\}$ labelled $\{1,...,6\}$ in lexicographic
order, we need only calculate $\Delta Z\{1,2\}$ and $\Delta Z\{1,6\}$.
The results are that $\Delta Z\{1,2\}\gg 0$ (has positive coefficients) and
that $\Delta Z\{1,6\} \gg (y_2y_5-y_3y_4)^2$.
In both cases $\Delta Z$ is nonnegative for all $\y>\zero$.
Thus $\II\mathsf{K}_4$ is Rayleigh.\\
Similar computations have shown that $\mathsf{K}_6$, $\mathsf{K}_{3,4}$
and $\mathsf{W_7}$ are $\II$--Rayleigh.\\
Alan Sokal showed me a computation that the Potts model
of $\mathsf{K}_4$ is Rayleigh for all $0<q\leq 1$, in Oct. 2005.
}\end{EG}

Conjecture 5.3 has been checked for $n\leq 6$.
\begin{CONJ}
If $e,f$ are adjacent edges in $\mathsf{K}_{n}$ then
$\Delta Z(\II \mathsf{K}_n)\{e,f\}$\\ $\gg 0$ (has positive coefficients).
\end{CONJ}

Conjecture 5.4 implies Conjecture 5.3.  The analogous statement
with $\BB \mathsf{G}$ in place of $\II \mathsf{G}$ follows from Theorem 5.6 of
Choe and Wagner \cite{CW}.
\begin{CONJ}
If $e,f$ are edges in a graph $\mathsf{G}$ and
$\Delta Z(\II \mathsf{G})\{e,f\}$ has negative coefficients, then
there is a $\mathsf{K}_{4}$--minor of $\mathsf{G}$ in which
$e$ and $f$ occur on non--adjacent edges.
\end{CONJ}
One might also conjecture that the only terms in $\Delta Z(\II\mathsf{G})\{e,f\}$
with negative coefficients already appear in $\Delta Z(\BB\mathsf{G})\{e,f\}$.
This is too optimistic -- it is true for $\mathsf{K}_4$
but not for $\mathsf{K}_5$.  Maybe this hints at a property of planar graphs,
but maybe not.\\

5.2. \ \textsc{two--sums of matroids.}\\

Assume that $\omega:\B(E)\goesto[0,\infty)$ and
$\nu:\B(F)\goesto[0,\infty)$ are Rayleigh. Then the
direct product $\omega\times\nu:\B(E\cup F)\goesto[0,\infty)$
is also Rayleigh.  This is a good exercise.

Let $M=Z(\M,q;\y)$ be the Potts model partition function of a
matroid $\M$, and let $g\in E= E(\M)$.
For any $S\subseteq E\drop\{g\}$,
$$\rank_{\M\drop g}(S)=\rank_{\M}(S),$$
from which it follows that
$$Z(\M\drop g,q;\y)=M|_{y_g=0}=M^g.$$
Also for any $S\subseteq E\drop\{g\}$,
$$\rank_{\M/g}(S)=\rank_{\M}(S\cup \{g\})-\rank_{\M}(\{g\}),$$
from which it follows that
$$Z(\M/g,q;\y)=q^{a(g)} \frac{\partial M}{\partial y_g},$$
in which $a(g)=0$ if $g$ is a loop of $\M$ and $a(g)=1$ otherwise.
For this reason it is convenient, for Potts models, to redefine
$M_g:=q^{a(g)} \partial M/\partial y_g$ in order that $M_g$ is the
Potts model of $\M/g$.

In the statement of Lemma 5.5, $\overline{S}$ denotes the closure of
$S\subseteq E(\M)$ in $\M$.
\begin{LMA}
Let $M=Z(\M,q;\y)$ be the Potts model partition function of a
matroid $\M$, and let $g\in E= E(\M)$ with $g$ not a loop of $\M$.\\
\textup{(a)}
$$M = M^{g}+q^{-1}y_{g}M_{g}.$$
\textup{(b)}
$$M^{g}-q^{-1}M_{g} =
(1-q^{-1})\sum_{S\subseteq E\drop\{g\}:\ g\not\in 
\overline{S}} q^{-\mathrm{rank}_{\M}(S)}\y^{S}.$$
\textup{(c)} If $q\neq 1$, then
$$\frac{M^{g}-M_{g}}{1-q} = \sum_{S\subseteq E\drop\{g\}:\ 
g\in\overline{S}} q^{-\mathrm{rank}_{\M}(S)}\y^{S}.$$
\textup{(d)} If $0<q<1$ and all $\y>\zero$, then
$q M^g < M_g \leq M^g$.  The weak inequality holds with equality if and
only if $g$ is a coloop of $\M$.
\end{LMA}
\begin{proof}
For part (a),
\showon
M
&=& \sum_{S\subseteq E:\ g\not\in S} q^{-\rank_\M(S)} \y^S
+ \sum_{S\subseteq E:\ g\in S} q^{-\rank_\M(S)} \y^S\\
&=&  \sum_{S\subseteq E\drop\{g\}} q^{-\rank_{\M}(S)} \y^S
+ y_g \sum_{S\subseteq E\drop\{g\}} q^{-\rank_{\M}(S\cup\{g\})} \y^S\\
&=& M^g + q^{-1} y_g M_g,
\showoff
as claimed.  For part (b),
\showon
M^{g}-q^{-1}M_{g}
&=&  \sum_{S\subseteq E\drop\{g\}}
\left[q^{-\rank_{\M}(S)}-q^{-\rank_{\M}(S\cup\{g\})}\right]\y^S\\
&=& (1-q^{-1})\sum_{S\subseteq E\drop\{g\}:\ g\not\in 
\overline{S}} q^{-\mathrm{rank}_{\M}(S)}\y^{S},\\
\showoff
as claimed.  For part (c),
\showon
\sum_{S\subseteq E\drop\{g\}:\ 
g\in\overline{S}} q^{-\mathrm{rank}_{\M}(S)}\y^{S}
&=&
M^g-\sum_{S\subseteq E\drop\{g\}:\ 
g\not\in\overline{S}} q^{-\mathrm{rank}_{\M}(S)}\y^{S}\\
&=& M^g - \frac{M^g - q^{-1} M_g}{1-q^{-1}}\\
&=& \frac{M^g - M_g}{1-q},
\showoff
as claimed.  For part (d), for every $S\subseteq E\drop\{g\}$,
$$\rank_\M(S)\leq \rank_\M(S\cup\{g\})\leq 1+\rank_\M(S).$$
Since $\rank_{\M/g}(S)=\rank_\M(S\cup\{g\})-1$ and $0<q<1$,
$$ q^{-\rank_{\M\drop g}(S)}\y^S\leq 
q^{-1-\rank_{\M/g}(S)}\y^S\leq
q^{-1-\rank_{\M\drop g}(S)}\y^S,$$
for every $S\subseteq E\drop\{g\}$.
Summing these inequalities shows that $q M^g \leq M_g \leq M^g$.
Since $g$ is not a loop of $\M$ the left inequality is strict.
The right inequality is tight if and only if
$\rank_\M(S\cup\{g\})= 1+\rank_\M(S)$ for every
$S\subseteq E\drop\{g\}$;\ that is, if and only if $g$ is a coloop of $\M$.
\end{proof}

The \emph{two--sum of matroids} $\N=\L\oplus_{g}\M$ is defined by means of
rank functions as follows.  Consider matroids $\L$ and $\M$ such that 
$E(\L)\cap E(\M)=\{g\}$, and such that $g$ is neither a loop nor a coloop in
 $\L$ or in $\M$.  Then $\N$ is the matroid on the set
$E(\N):=E(\L)\cup E(\M)\drop\{g\}$ with rank function
$$\rank_\N (S) := \rank_\L (S\cap E(\L)) + \rank_\M (S\cap E(\M)) - \nu(S)$$
for all $S\subseteq E(\N)$, in which
$$
\nu(S):=
\left\{\begin{array}{ll}
1 & \mathrm{if}\  g\in \overline{S\cap E(\L)}\ \mathrm{and}\ 
    g\in \overline{S\cap E(\M)},\\
0 & \mathrm{otherwise}.
\end{array}\right.
$$

\begin{PROP}
Consider the two--sum of matroids $\N=\L\oplus_{g}\M$.
Write $L=Z(\L,q;\y)$, $M=Z(\M,q;\y)$ and $N=Z(\N,q;\y)$
for the corresponding Potts model partition functions.\\
\showon
N &=& 
L^{g}M^{g}-\frac{1}{1-q}\left(L^{g}-L_{g}\right)\left(M^{g}-M_{g}\right).\\
N &=&
\frac{1}{1-q}(-qL^{g}M^{g}+L^{g}M_{g}+L_{g}M^{g}-L_{g}M_{g}).
\showoff
\end{PROP}
\begin{proof}
From Lemma 5.5 it follows that
\showon
N
&=& \sum_{S\subseteq E(\N)} q^{-\rank_\N(S)}\y^S\\
&=& \sum_{S\subseteq E(\L)\drop\{g\}}\,\,
    \sum_{T\subseteq E(\M)\drop\{g\}}
	q^{-\rank_\L(S)-\rank_\M(T)+\nu(S\cup T)}\y^{S\cup T}\\
&=& L^g M^g -(1-q)\frac{L^g-L_g}{1-q}\cdot \frac{M^g-M_g}{1-q}.
\showoff
The second equation follows by routine algebra.
\end{proof}

The limiting argument of Example 3.2 implies the following.
\begin{CORO}
Consider the two--sum of matroids $\N=\L\oplus_{g}\M$.\\
\textup{(a)}\ For bases, $N=Z(\BB\N;\y)$ et cetera,
$$N=L^gM_g + L_g M^g.$$
\textup{(b)}\ For independent sets, $N=Z(\II\N;\y)$ et cetera,
$$N=L^g M_g + L_g M^g - L_g M_g.$$
\textup{(c)}\ For spanning sets, $N=Z(\SS\N;\y)$ et cetera,
$$N=L^g M_g + L_g M^g - L^g M^g.$$
\end{CORO}

The Potts--Rayleigh condition is that $M=Z(\M,q;\y)$ is Rayleigh
for all $q$ in some interval $0<q\leq q_*(\M)$.
Is it true that if $M$ is Rayleigh at $q=q_0$ then $M$ is
Rayleigh for all $0<q\leq q_0$?
Define $q_c(\M)$ to be the supremum of all $0\leq q< 1$ for which
$Z(\M,q;\y)$ is Rayleigh. 
The argument of Lemma 4.1 can be adapted to show that
for all $g\in E(\M)$, $q_c(\M)\leq q_c(\M/g)$ and
$q_c(\M)\leq q_c(\M\drop g)$.

\begin{THM}
The following classes of matroids are closed by taking two--sums:\\
\textup{(a)}\ $\BB$--Rayleigh matroids;\\
\textup{(b)}\ $\II$--Rayleigh matroids;\\
\textup{(c)}\ $\SS$--Rayleigh matroids;\\
\textup{(d)}\ Potts--Rayleigh matroids.
\end{THM}
\begin{proof}
Part (a) is Theorem 3.5 of Choe and Wagner \cite{CW}.

Parts (b) and (c) are equivalent by $\M\leftrightarrow\M^*$ duality.
To prove (b) we repeat the argument for part (a) in a more complicated setting.
Consider a two--sum of matroids $\N=\L\oplus_g\M$
in which $\II\L$ and $\II\M$ are Rayleigh, and let $N=Z(\II\N;\y)$,
\emph{et cetera}.  From Corollary 5.7(b),
$N=L^g M_g + L_g M^g - L_g M_g.$  Fix two distinct
elements $e\neq f$ in $E(\N)$, and $\y>\zero$.  There are two cases
(by symmetry) for $\Delta N\{e,f\}$:\
either $\{e,f\}\subseteq E(\L)\drop\{g\}$,
or $e\in E(\L)\drop\{g\}$ and $f\in E(\M)\drop\{g\}$.

If $e\in E(\L)\drop\{g\}$ and $f\in E(\M)\drop\{g\}$ then
a moderately taxing calculation yields
$$\Delta N\{e,f\}= \Delta L\{e,g\} \cdot \Delta M\{g,f\}.$$
(We skip the details since an analogous calculation occurs in the
proof of part (d).)  This is nonnegative since both $L$ and $M$ are Rayleigh.

If $\{e,f\}\subseteq E(\L^g)$ then
\showon
& & \Delta N\{e,f\}(\y)\\
&=& (M^g)^2 \left[ \Delta L^g\{e,f\}(\y) + y_g\Theta L\{e,f|g\}(\y)
+y_g^2 \Delta L_g\{e,f\}(\y) \right] \\
&=& (M^g)^2 \Delta L\{e,f\}(\y)
\showoff
by setting $y_g=M_g/M^g>0$.
(Again we skip the details since an analogous but much more difficult
calculation occurs in the proof of part (d).)
Since $M^g(\y)>0$, this is well--defined.
Since $\II\L$ is Rayleigh, this is nonnegative.
It follows that $\II\N$ is Rayleigh.

For part (d) we repeat the argument a third time.
Consider a two--sum of matroids $\N=\L\oplus_g\M$
for which the Potts models $L=Z(\L,q;\y)$ and $M=Z(\M,q;\y)$
are Rayleigh on intervals $0<q\leq q_*(\L)$ and $0<q\leq q_*(\M)$
respectively.  We use the formula
$$ N =
\frac{1}{1-q}(-qL^{g}M^{g}+L^{g}M_{g}+L_{g}M^{g}-L_{g}M_{g})
$$
of Proposition 5.6.

Fix two distinct elements $e\neq f$ in $E(\N)$, and $\y>\zero$.  There are
two cases (by symmetry) for $\Delta N\{e,f\}$:\
either $\{e,f\}\subseteq E(\L)\drop\{g\}$,
or $e\in E(\L)\drop\{g\}$ and $f\in E(\M)\drop\{g\}$.

If $e\in E(\L)\drop\{g\}$ and $f\in E(\M)\drop\{g\}$ then
\showon
(1-q)^2 N_e^f N_f^e
&=&  
(-qL_e^{g}M^{fg}+L_e^{g}M_{g}^f+L_{eg}M^{fg}-L_{eg}M_{g}^f)\\
& & \cdot(-qL^{eg}M_f^{g}+L^{eg}M_{fg}+L_{g}^eM_f^{g}-L_{g}^e M_{fg})
\showoff
and
\showon
(1-q)^2 N_{ef} N^{ef}
&=&
(-qL_e^{g}M_f^g+L_e^{g}M_{fg}+L_{eg}M_f^g-L_{eg}M_{fg})\\
& & \cdot(-qL^{eg}M^{fg}+L^{eg}M_{g}^f+L_{g}^e M^{fg}-L_{g}^e M_{g}^f).
\showoff
The amazing fact is that these $32$ terms cancel almost completely in the
difference $(1-q)^2\Delta N\{e,f\}$.  The cancellation is exact except for
terms with two $g$s up and two $g$s down in the deletion/contraction notation.
Moreover, the four remaining middle terms factor as
$$\Delta N\{e,f\} = \frac{1}{1-q}\Delta L\{e,g\}\cdot\Delta M\{g,f\}.$$
For $0<q\leq \min\{q_*(\L),q_*(\M)\}$
this is nonnegative since both $L$ and $M$ are Rayleigh for $q$ in this range.

If $\{e,f\}\subseteq E(\L)\drop\{g\}$ then
\showon
(1-q)^2 N_e^f N_f^e
&=&  ( -q L_e^{fg} M^g + L_e^{fg} M_g + L_{eg}^f M^g -L_{eg}^f M_g)\\
& & \cdot  ( -q L_f^{eg} M^g + L_f^{eg} M_g + L_{fg}^e M^g -L_{fg}^e M_g)
\showoff
and
\showon
(1-q)^2 N_{ef} N^{ef}
&=&  ( -q L_{ef}^g M^g + L_{ef}^g M_g + L_{efg} M^g -L_{efg} M_g)\\
& & \cdot  ( -q L^{efg} M^g + L^{efg} M_g + L_g^{ef} M^g -L_g^{ef} M_g).
\showoff

Let's collect $(1-q)^2\Delta N\{e,f\}$ according to powers of $M_g/M^g$:
\showon
& & (1-q)^2 \Delta N\{e,f\}=\\
& &(M^g)^2
\left[  (L_{eg}^f-q L_e^{fg})(L_{fg}^e-q L_f^{eg})
-        (L_{efg}-q L_{ef}^{g}) (L_{g}^{ef}-q L^{efg}) \right]\\
&+&(M_g M^g) 
\left[  (L_{eg}^f-q L_e^{fg}) (L_f^{eg}-L_{fg}^e)
 + (L_e^{fg}-L_{eg}^f) (L_{fg}^e-q L_f^{eg})\right]\\
&-&(M_g M^g)
\left[ (L_{efg}-q L_{ef}^{g}) (L^{efg}-L_{g}^{ef})
 +  (L_{ef}^{g}-L_{efg}) (L_{g}^{ef}-q L^{efg})\right]\\
&+&(M_g)^2 
\left[  (L_e^{fg}-L_{eg}^f)(L_f^{eg}-L_{fg}^e)
-(L_{ef}^{g}-L_{efg})(L^{efg}-L^{ef}_{g}) \right].
\showoff
Let $\gamma:=M_g/M^g$.  This is a nonzero rational function of $q$ and
$\{y_c:\ c\in E(\M\drop g)\}$ with positive coefficients.
When $0<q<1$ and $\y>\zero$ it follows from Lemma 5.5(d)
that $q<\gamma<1$, since $g$ is not a loop or coloop in $\L$.
Simplifying the above formula we obtain
$$ (1-q)^2\Delta N\{e,f\} = (M^g)^2(C+B \gamma +A \gamma^2) $$
in which
\showon
C &=&
q^2\Delta L^g\{e,f\} -q \Theta L\{e,f|g\} + \Delta L_g\{e,f\}   \\
B &=&
-2q \Delta L^g\{e,f\} +(1+q)\Theta L\{e,f|g\} -2 \Delta L_g\{e,f\}\\
A &=&
\Delta L^g\{e,f\} - \Theta L\{e,f|g\} + \Delta L_g\{e,f\}.
\showoff
Another way to write this is 
\showon
& & C+B \gamma+A \gamma^2\\
&=& (\gamma-q)^2 \left[ \Delta L^g\{e,f\} 
+y_g \Theta L\{e,f|g\}
+y_g^2 \Delta L_g\{e,f\} \right].
\showoff
in which $y_g:=(1-\gamma)/(\gamma-q)$ is positive and finite for all
$0<q< 1$ and $\y>\zero$.  Thus
$$(1-q)^2\Delta N\{e,f\}=
(M^g)^2(\gamma-q)^2 \Delta L\{e,f\}(y_g=(1-\gamma)/(\gamma-q)).$$
This is nonnegative for all
$0<q\leq q_*(\L)$.
Therefore, the Potts model of $\N=\L\oplus_g\M$ is Rayleigh for
$$0<q\leq \min\{q_*(\L), q_*(\M)\}.$$
This completes the proof.
\end{proof}

\vspace{4mm}

5.3. \ \textsc{classes of rayleigh matroids.}\\

\begin{PROP}
Let $\mathsf{G}=(V,E)$ be a finite connected graph.  For each
$S\in\II\mathsf{G}$ let $\omega_G(S)$ be the product of the sizes of the
connected components of $(V,S)$.  Then $\omega_G$, supported on $\II\mathsf{G}$,
is such that $\sum_{k=0}^{n-1} f_k(\omega_G) t^{n-k}$ has only real (nonpositive)
zeros, the condition \textup{(f-6)}.
\end{PROP}
\begin{proof}
Orient the edges of $\mathsf{G}$ arbitrarily, and let $D$ be the 
corresponding $V$--by--$E$ signed incidence matrix, and $Y:=\diag(y_e:\
e\in E)$, and $Q:=DYD^\dagger$. 
For any $\y>\zero$ the matrix $Q$ is Hermitian, in
fact positive semidefinite with a one--dimensional nullspace,
so that $\det(tI+Q)$ has only real zeros.
By the Binet--Cauchy identity and the principal minors matrix--tree theorem
\cite{Cha},
\showon
\det(tI+Q)
&=& \sum_{R\subseteq V} t^{|R|} \det Q(R|R)\\
&=& \sum_{R\subseteq V} t^{|R|}\sum_{S\in\B(E,n-|R|)} |\det D(R|S]|^2\y^S \\
&=& \sum_{k=0}^n t^{n-k} \sum_{S\in\B(E,k)} \omega_G(S)\y^S.
\showoff
For all $\y>\zero$ this polynomial has only real (nonpositive) zeros.
When $\y=\one$ this implies (f-6).
\end{proof}

\begin{CONJ}
For any finite connected graph $\mathsf{G}$,
the $\omega_G$ defined in Proposition 5.9 is Rayleigh.
\end{CONJ}

Here is a very interesting formula for $Z(\II\mathsf{G};\y)$
from which we might be able to see that $\II\mathsf{G}$ is Rayleigh.
The Grassmann--Berezin calculus technique of \cite{Ab,BI,CJSSS} is
great, and in particular formula (13) of CJSSS \cite{CJSSS} implies that
$$t^nZ(\II\mathsf{G};\y/t)=\int \mathrm{d}(\psi\bar{\psi})\,
\exp\left[\bar{\psi}(tI+Q)\psi -
t \sum_{e=\la ij \ra \in E} y_e \bar{\psi}_i\psi_i \bar{\psi}_j\psi_j
\right].$$
Here $\psi_i$ and $\bar{\psi}_i$ are fermionic degrees of freedom
associated with vertex $i$,
and $Q=DYD^\dagger$ is as in Proposition 5.9.
Combinatorially, the idea is that trees can also be rooted negatively at their
edges as well as positively at their vertices. Since a tree has one more vertex
than edge, each tree gets a net count of one, and so each spanning
forest is counted exactly once.
It should be possible to compute $\Delta Z(\II\mathsf{G})\{e,f\}$ from this --
to see if it is positive will be more difficult. 
Generalization of this formula to wider classes of matroids is also a worthy goal.
For comparison, note that
$$\det(tI+Q)= \int \mathrm{d}(\psi\bar{\psi})\,
\exp\left[\bar{\psi}(tI+Q)\psi\right].$$

To conclude, let's sum up our paltry stock of examples.

The class of Potts--Rayleigh matroids:\\
* contains all uniform matroids $\U(m,r)$;\\
* contains the graph $\mathsf{K}_4$;\\
* is closed by taking duals, minors, and two--sums.

The class of $\II$--Rayleigh matroids:\\
* contains all Potts--Rayleigh matroids;\\
* contains the graphs $\mathsf{K}_{6}$, $\mathsf{K}_{3,4}$, and $\mathsf{W}_7$;\\
* is closed by taking minors and two--sums.

The class of Potts--Rayleigh matroids contains all series--parallel graphs.
Thus, if the Big Conjecture $3.4$ is true then every series--parallel graph
satisfies Mason's Conjecture (I-4).

\begin{CONJ}
The class of $\II$--Rayleigh matroids contains\\
\textup{1.}\ all planar graphs;\\
\textup{2.}\ all graphs;\\
\textup{3.}\ all regular matroids;\\
\textup{4.}\ all sixth--root of unity matroids;\\
\textup{5.}\ all half--plane property matroids;\\
\textup{6.}\ all $\BB$--Rayleigh matroids.
\end{CONJ}
By familiar results and others in \cite{COSW,CW}, these conjectures
are increasingly strong.
Conjecture 5.11.6 is that $\II\M$ is Rayleigh if and only if $\BB\M$ is
Rayleigh -- this is almost certainly false, and a specific counterexample
is much to be desired.  Conjecture 5.11.2 seems reasonable.
Similar conjectures could be made about the Potts--Rayleigh
condition, with even less evidence.

\end{document}